\documentclass[12pt,a4paper]{amsart}

\usepackage[utf8]{inputenc} % allow utf-8 input
\usepackage[T1]{fontenc}    % use 8-bit T1 fonts

\usepackage{epic}
\numberwithin{equation}{section}
\usepackage[margin=2.9cm]{geometry}

\usepackage{amsaddr,amsmath,amssymb,amsfonts,array,dsfont,amsthm}
\usepackage{graphicx,subcaption}

\newtheorem{theorem}{Theorem}
\newtheorem{lemma}{Lemma}

\theoremstyle{plain}
\newtheorem*{graph_minor_theorem}{Graph Minor Theorem}
\newtheorem*{kuratowski_theorem}{Kuratowski's Theorem}
\newtheorem*{wagner_theorem}{Wagner's Theorem}

\title{Non-separating Planar Graphs}

\author{Hooman R. Dehkordi}
\address{Faculty of Information Technology, Monash University, \\
Clayton, Victoria 3800, Australia \\
\texttt{hooman.dehkordi@monash.edu}}

\author{Graham Farr}
\address{Faculty of Information Technology, Monash University, \\
Clayton, Victoria 3800, Australia \\
\texttt{graham.farr@monash.edu}}

\begin{document}
\maketitle

\begin{abstract}
A graph $G$ is a non-separating planar graph if there is a drawing $D$ of $G$ on the plane such that (1) no two edges cross each other in $D$ and (2) for any cycle $C$ in $D$, any two vertices not in $C$ are on the same side of $C$ in $D$. 

Non-separating planar graphs are closed under taking minors and are a subclass of planar graphs and a superclass of outerplanar graphs.

In this paper, we show that a graph is a non-separating planar graph if and only if it does not contain
$K_1 \cup K_4$ or $K_1 \cup K_{2,3}$ or $K_{1,1,3}$ as a minor.

Furthermore, we provide a structural characterisation of this class of graphs. More specifically, we show that any maximal non-separating planar graph is either an outerplanar graph or a subgraph of a wheel or it can be obtained by subdividing some of the side-edges of the 1-skeleton of a triangular prism (two disjoint triangles linked by a perfect matching).

Lastly, to demonstrate an application of non-separating planar graphs, we use the characterisation of non-separating planar graphs to prove that there are maximal linkless graphs with $3n-3$ edges which provides an answer to a question asked by Horst Sachs about the number of edges of linkless graphs in 1983. 
%The previous best known upperbound on the number of edges of maximal linkless graphs was $4n-10$ edges.
\end{abstract}

% keywords can be removed
%\keywords{First keyword \and Second keyword \and More}

\section{Introduction}\label{se:intro}

A \emph{drawing} of a graph in the plane consists of a set of points representing (under a bijection) the vertices of the graph and a set of curves between certain pairs of points representing edges between corresponding vertex pairs of the graph where the curves do not pass through the points that represent vertices. A \emph{planar drawing} is a drawing in which edges do not  intersect.

Let $C$ be a cycle in a planar drawing $D$ of a graph $G$, then $C$ is a \emph{separating cycle} if there is at least one vertex in the interior of $C$ and one vertex in the exterior of $C$.

\begin{figure}[h]
  \centering
  \begin{subfigure}[b]{0.27\textwidth}
%    \hspace{7 mm} % pushing figure to the right so that it is above the sub caption
    \includegraphics[width=0.80\textwidth]{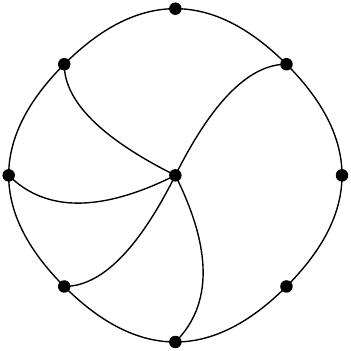} % scaling the fig so that the size of vertices is not too big
%    \caption{$K_1 \cup K_4$}
    \label{sfig:K_1_and_K_4}
  \end{subfigure}
  \hspace{5 mm}
  \begin{subfigure}[t]{0.27\textwidth}
%    \hspace{7 mm} % pushing figure to the right so that it is above the sub caption
    \includegraphics[width=0.9\textwidth]{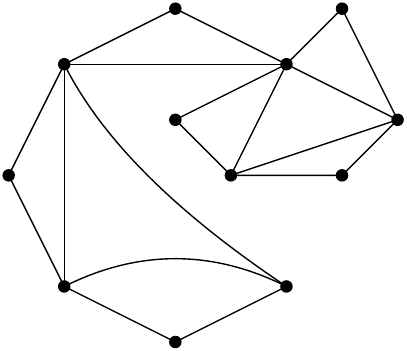} % scaling the fig so that the size of vertices is not too big
%    \caption{$K_1 \cup K_{2,3}$}
    \label{sfig:K_1_and_K_23}
  \end{subfigure}
  \hspace{5 mm}
  \begin{subfigure}[t]{0.27\textwidth}
    \hspace{7 mm} % pushing figure 3mm to the right so that it is above the sub caption
    \includegraphics[width=0.74\textwidth]{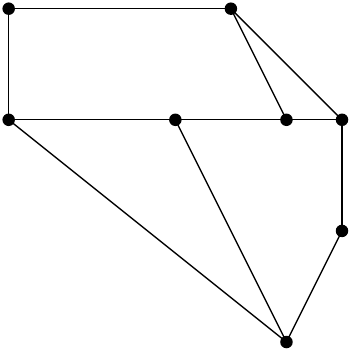}
%    \caption{$K_{1,1,3}$}
    \label{sfig:K_113}
  \end{subfigure}
  \caption{Three examples of non-separating planar graphs}
  \label{fig:non-separating_planar_graphs_example}
\end{figure}

A \emph{non-separating planar drawing} of a graph is a planar drawing of the graph that does not contain any separating cycles. A \emph{non-separating planar graph} is a graph that has a non-separating planar drawing (see for example Figure~\ref{fig:non-separating_planar_graphs_example}).

Any graph $G'$ that can be obtained from a graph $G$ by a series of edge deletions, vertex deletions and edge contractions is called a \emph{minor} of $G$.

A set $S$ of graphs is a \emph{minor-closed set} or \emph{minor-closed family} of graphs if any minor of a graph $G \in S$ is also a member of $S$. 

In this paper we characterise non-separating planar graphs. Non-separating planar graphs are a subclass of planar graphs and a superclass of outerplanar graphs and are closed under minors. To characterise non-separating planar graphs we prove Theorems~\ref{th:non-separating_planar_graphs_excluded_minors} and \ref{th:structural_characterisation_of_non-separating_planar_graphs} as follows.

\begin{figure}[h]
  \centering
  \begin{subfigure}[t]{0.27\textwidth}
%    \hspace{7 mm} % pushing figure to the right so that it is above the sub caption
    \includegraphics[width=0.80\textwidth]{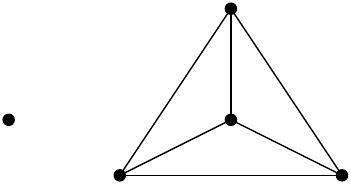} % scaling the fig so that the size of vertices is not too big
    \caption{$K_1 \cup K_4$}
    \label{sfig:K_1_and_K_4}
  \end{subfigure}
  \hspace{5 mm}
  \begin{subfigure}[t]{0.27\textwidth}
%    \hspace{7 mm} % pushing figure to the right so that it is above the sub caption
    \includegraphics[width=0.80\textwidth]{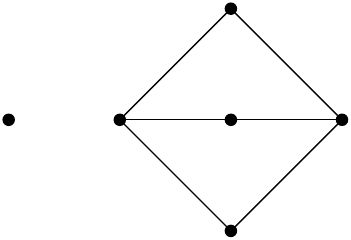} % scaling the fig so that the size of vertices is not too big
    \caption{$K_1 \cup K_{2,3}$}
    \label{sfig:K_1_and_K_23}
  \end{subfigure}
  \hspace{5 mm}
  \begin{subfigure}[t]{0.27\textwidth}
    \hspace{7 mm} % pushing figure 3mm to the right so that it is above the sub caption
    \includegraphics[width=0.57\textwidth]{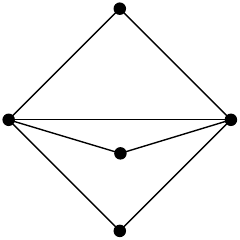}
    \caption{$K_{1,1,3}$}
    \label{sfig:K_113}
  \end{subfigure}
  \caption{Excluded minors for non-separating planar graphs}
  \label{fig:excluded_minors_for_non_separating_planar_graphs}
\end{figure}

\begin{theorem}
\label{th:non-separating_planar_graphs_excluded_minors}
	A graph $G$ is a non-separating planar graph if and only if it does not contain any of $K_1 \cup K_4$ or $K_1 \cup K_{2,3}$ or $K_{1,1,3}$ as a minor\footnote{where $\cup$ denotes the disjoint union} (see Figure~\ref{fig:excluded_minors_for_non_separating_planar_graphs}). 
\end{theorem}

\begin{figure}
  \centering
  \begin{subfigure}[t]{0.40\textwidth}
%    \hspace{7 mm} % pushing figure to the right so that it is above the sub caption
    \includegraphics[width=0.7\textwidth]{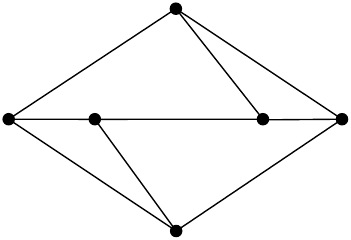} % scaling the fig so that the size of vertices is not too big
    \caption{Triangular prism}
    \label{sfig:triangular_prism}
  \end{subfigure}
  \hspace{2 mm}
  \begin{subfigure}[t]{0.40\textwidth}
    \hspace{7 mm} % pushing figure 3mm to the right so that it is above the sub caption
    \includegraphics[width=0.70\textwidth]{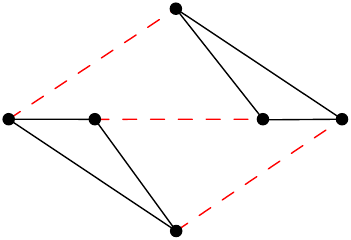}
    \caption{Elongated triangular prism}
    \label{sfig:elongated_triangular_prism}
  \end{subfigure}
  \caption{Triangular prism and elongated triangular prism}
  \label{fig:triangular_prism_and_elongated_triangular_prism}
\end{figure}

An edge $e = (u, v)$ in a graph $G$ is \emph{subdivided} by replacing it with two edges $(u,w), (w,v)$ where $w$ is not a vertex of $G$. A \emph{subdivision} of a graph $G$ is a graph that can be obtained by some sequence of subdivisions, starting with $G$. A graph is a \emph{triangular prism} if it is isomorphic to the graph that is depicted in Figure~\ref{sfig:triangular_prism}. A graph is an \emph{elongated triangular prism} if it is a triangular prism or if it is obtained by some sequence of subdivisions of the red dashed edges of the triangular prism depicted in Figure~\ref{sfig:elongated_triangular_prism}.

An \emph{outerplanar drawing} is a drawing of a graph on a disk in which no two edges cross and and all the vertices of the graph are located on the boundary of the disk. \emph{Outerplanar graphs} are the graphs that have an outerplanar drawing. 

We also characterise non-separating planar graphs in terms of their structure as follows.

\begin{theorem}
\label{th:structural_characterisation_of_non-separating_planar_graphs}
	Any non-separating planar graph is one of the following:
	\begin{enumerate}
		\item an outerplanar graph,
		\item a subgraph of a wheel,
		\item a subgraph of an elongated triangular prism.
	\end{enumerate}
\end{theorem}

A \emph{realisation} $\mathcal{R}$ of a graph $G = (V, E)$ in $\mathds{R}^3$ consists of a set of points in $\mathds{R}^3$ that represent (under a bijection) the vertices of the graph and a set of curves between certain pairs of points that represent the edges between corresponding vertex pairs of the graph, such that these curves do not intersect and also do not pass through the points that represent the vertices of the graph. Informally, realisations of graphs are drawings of graphs in $\mathds{R}^3$.

\begin{figure}[tbp]
  \centering
  \includegraphics[width=3cm]{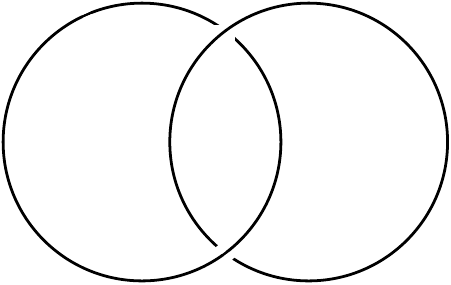}
  \caption{A link is composed of two cycles in three dimensions that cannot be separated from each other.}
  \label{fig:link}
\end{figure}

Two vertex-disjoint cycles $C_1$ and $C_2$ that are embedded into $\mathds{R}^3$ are \emph{linked} if no topological sphere can be embedded into $\mathds{R}^3$ separating $C_1$ from $C_2$. Two linked cycles are called a \emph{link}. To put it in another way, two cycles $C_1$ and $C_2$ are not linked (unlinked) if they can be continuously deformed without ever intersecting each other until $C_1$ and $C_2$ end up in two different sides of a topological sphere embedded into $\mathds{R}^3$. Informally, a link consists of two cycles that are embedded in three dimensions such that they cannot be separated unless we cut one of them (see Figure~\ref{fig:link}). 

A realisation $\mathcal{R}$ of a graph is \emph{linkless} if it contains no links. A graph is \emph{linkless} if it has a linkless realisation. Although linkless graphs are characterised in terms of a set of forbidden minors, there are a lot of unanswered questions about them. For example, since linkless graphs do not contain a $K_6$-minor, it follows that they have at most $4|V|-10$ edges where $|V|$ is the number of vertices in $G$ \cite{mad68}. However, we do not know whether all maximal linkless graphs have $4|V|-10$ edges. In this paper we show that there is a class $\mathcal{G}$ of maximal linkless graphs such that any graph $G \in \mathcal{G}$ has at most $3|V| - 3$ edges which provides an answer to a question asked by Horst Sachs about the number of edges of linkless graphs in 1983~\cite{sach83}. More specifically, we prove the following theorem:

\begin{theorem}
\label{th:number_of_edges_in_flat_graphs}
	There exists an infinite family $\mathcal{G}$ of maximal linkless graphs such that any graph $G \in \mathcal{G}$ has at most $3|V(G)| - 3$ edges.
\end{theorem}

The rest of this paper is organised as follows. Section~\ref{se:background} goes into more details about the different classes of graphs that we are dealing with in this paper and describes their relation to each other. Section~\ref{se:lemmas} is dedicated to proving a number of preliminary lemmas that are later used in proofs of the main theorems of this paper. More specifically, in Section~\ref{se:lemmas} we investigate the structure of the graphs that contain $K_{2,3}$ as a minor but do not contain any of $K_1 \cup K_{2,3}$, $K_1 \cup K_4$ or $K_{1,1,3}$ as a minor. These lemmas are structured in this specific manner to be useful in characterising non-separating planar graphs both in terms of forbidden minors and also in terms of their structure. Section~\ref{se:proof} then uses the results of Section~\ref{se:lemmas} to prove the main theorems of this paper. In Section ~\ref{se:linkless_edges_count} we demonstrate the relationship of non-separating planar graphs with linkless graphs by using Theorem~\ref{th:non-separating_planar_graphs_excluded_minors} to prove Theorem~\ref{th:number_of_edges_in_flat_graphs}. Lastly, in Section~\ref{se:conclusion}, we summarise our results and point out future directions for research.

%%%%%%%%%%%%%%%%%%%%%%%%%%%%%%%%%%%%%%%%
%
%  Background
%
%%%%%%%%%%%%%%%%%%%%%%%%%%%%%%%%%%%%%%%%

\section{Background}
\label{se:background}

The theory of graph minors developed by Robertson and Seymour is one of the most important recent advances in graph theory and combinatorics. This substantial body of work is presented in a series of 23 papers (Graph Minors I--XXIII) over 20 years from 1983 to 2004.

The graph minor theorem (also known as Robertson--Seymour theorem or Wagner's conjecture) can be formulated as follows:
\begin{graph_minor_theorem} [Robertson and Seymour \cite{rob04}]
\label{co:graph_minor_theorem_1}
	Every minor-closed class of graphs can be characterised by a finite family of excluded minors.
\end{graph_minor_theorem}

Perhaps the most famous minor-closed class of graphs is the class of planar graphs. \emph{Planar graphs} are the graphs that have a planar drawing.

In 1930, Kuratowski characterised planar graphs in terms of two forbidden subdivisions.
More specifically, he proved the following theorem: 
\begin{kuratowski_theorem}[Kuratowki \cite{kur30}]
    A graph is planar if and only if it does not contain a subdivision of $K_5$ or a subdivision of $K_{3,3}$ as a subgraph.
    \label{th:kuratowski}
\end{kuratowski_theorem}
%This theorem played an important role in the development of two-dimensional graph drawing. 
%These two forbidden graphs are depicted in Figure \ref{fig:forbidden_minors_of_planar_graphs}. 
Later on, Wagner characterised planar graphs in terms of forbidden minors as follows: 
\begin{wagner_theorem}[Wagner \cite{wag36}]
    A graph is planar if and only if it does not have $K_5$ or $K_{3,3}$ as a minor.
    \label{th:wagner}
\end{wagner_theorem}

In fact, it is easy to see that for any surface\footnote{a 2-manifold, or in other words a topological space such that every point has a neighborhood that is homeomorphic to an open subset of a Euclidean plane} $\Sigma$, the class of graphs that can be drawn on $\Sigma$ without edge crossings is closed under minors. For example, the class of toroidal graphs (graphs that can be drawn on a torus without edge crossings) is also closed under minors and hence can be characterised in terms of a finite set of forbidden minors. However, the complete set of forbidden minors for this class of graphs is not yet known~\cite{myrvold2018large}.

Another wellknown minor-closed class of graphs is the class of outerplanar graphs. Chartrand and Harary proved that a graph is outerplanar if and only if it does not contain $K_4$ or $K_{2,3}$ as a minor \cite{char67}.

%A \emph{disked realisation}\index{disked realisation}\marginnotes{disked realisation} $\mathcal{R}$ of a graph $G$ is a realisation of $G$, in which for each cycle $C$ in $\mathcal{R}$, there is an open disk $d$ in $\mathcal{R}$ that is bounded by $C$. A disked realisation $\mathcal{R}$ of a graph $G$ is \emph{flat}\index{flat realisation}\marginnotes{flat realisation} if every disk $d$ in $\mathcal{R}$ is disjoint from the vertices and edges of $G$. A \emph{flat graph}\index{flat graph}\marginnotes{flat graph} is a flat embedding of a graph. \todo{Fix the definition of flat graphs. In fact this is the embedding of the graph of a flat disk realisation.}
%

%A realisation $\mathcal{R}$ of a graph $G$ is \emph{flat} if for every cycle $C$ in $\mathcal{R}$ we can embed an open disk $d$ such that $d$ is bounded by $C$ and is disjoint from the vertices and edges of $G$. A \emph{flat graph} is a flat embedding of a graph.

%It is easy to see that any flat graph is a linkless graph. Robertson, Seymour and Thomas proved that any linkless graph has a flat embedding and therefore the class of linkless and flat graphs are equivalent \cite{rob91,rob95}. Moreover, 

Linkless graphs are also closed under minors. Sachs suggested the study of linkless embeddings for the first time \cite{sach83}. He conjectured that these embeddings can be characterised by excluding the Petersen family of graphs.

\begin{figure}
  \centering
  \begin{subfigure}[t]{0.24\textwidth}
%    \hspace{7 mm} % pushing figure 3mm to the right so that it is above the sub caption
    \includegraphics[width=\textwidth]{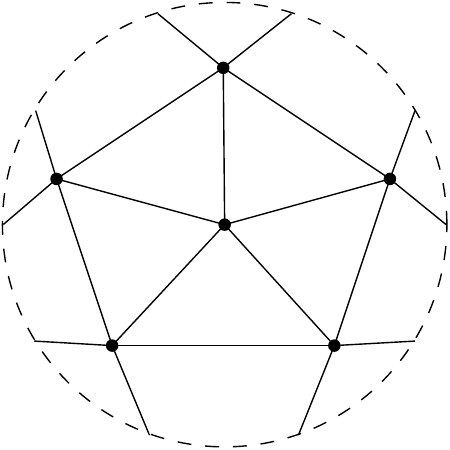} % scaling the fig so that the size of vertices is not too big
    \caption{$K_6$}
    \label{sfig:0_p(K_6)}
  \end{subfigure}
%  \hspace{3 mm}
  \begin{subfigure}[t]{0.24\textwidth}
    \includegraphics[width=\textwidth]{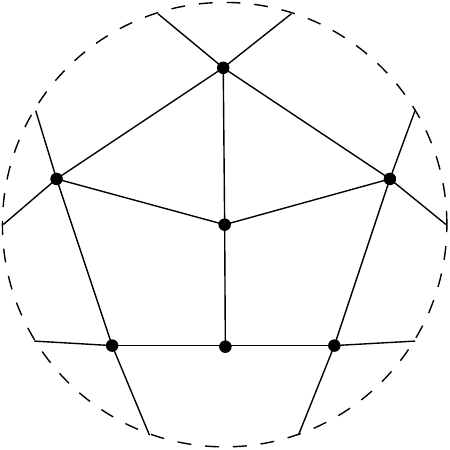}
    \caption{}
    \label{sfig:1_p}
  \end{subfigure}
%  \hspace{3 mm}
  \begin{subfigure}[t]{0.25\textwidth}
    \includegraphics[width=\textwidth]{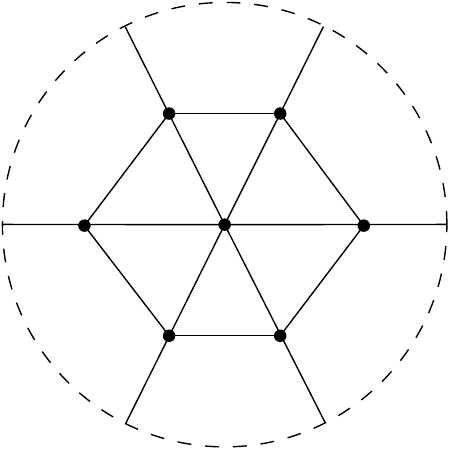}
    \caption{}
    \label{sfig:3_p(K_1,3,3)}
  \end{subfigure}
%  \hspace{3 mm}
  \begin{subfigure}[t]{0.24\textwidth}
    \includegraphics[width=\textwidth]{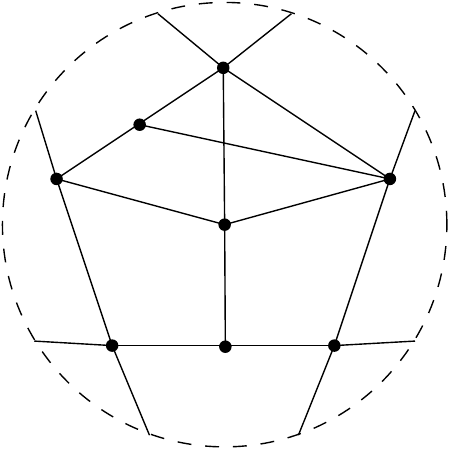}
    \caption{}
    \label{sfig:3_3(K_1,3,3)}
  \end{subfigure}
  \\
%  \hspace{3 mm}
  \begin{subfigure}[t]{0.25\textwidth}
    \includegraphics[width=\textwidth]{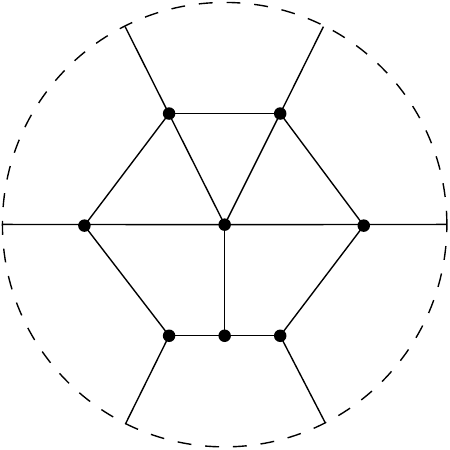}
    \caption{}
    \label{sfig:2_p}
  \end{subfigure}
%  \hspace{3 mm}
  \begin{subfigure}[t]{0.25\textwidth}
    \includegraphics[width=\textwidth]{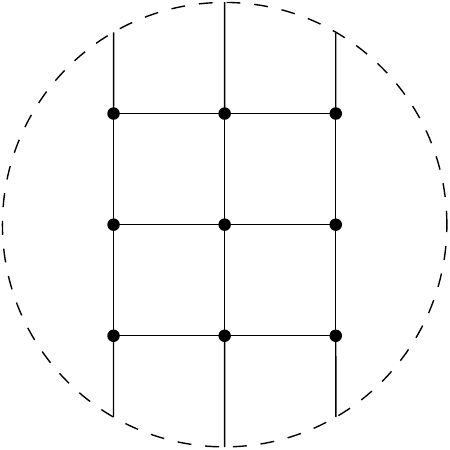}
    \caption{}
    \label{sfig:6_p}
  \end{subfigure}
%  \hspace{3 mm}
  \begin{subfigure}[t]{0.25\textwidth}
    \includegraphics[width=\textwidth]{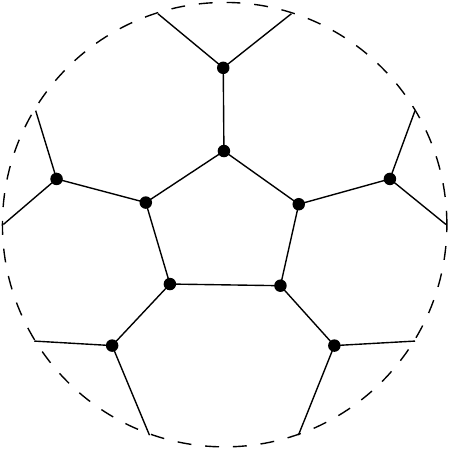}
    \caption{Petersen graph}
    \label{sfig:5_p}
  \end{subfigure}
  \caption{The Petersen family of graphs on the projective plane. (These drawings are drawn based on the drawings in \cite{kaw10}.)}
  \label{fig:Petersen_family}
\end{figure}

The \emph{Petersen family} of graphs consists of $K_6$ and six other graphs including the Petersen graph, as shown in Figure \ref{fig:Petersen_family}.

Conwoy, Gordon and Sachs proved that $K_6$ is not linkless \cite{con83}. Sachs has also proved that the other members of Petersen family of graphs are not linkless \cite{sach83}. Moreover in the same paper he showed that every minor of a linkless graph is linkless.

Robertson, Seymour and Thomas proved that $G$ is linklessly embeddable in $\mathds{R}^3$ if and only it does not contain any graph in the Petersen family as a minor \cite{rob91,rob95}. 

Among the other characterisations of graphs in terms of forbidden minors, we point out the following famous results: 
\begin{itemize}
	\item characterisation of the projective planar graphs (graphs that are embeddable on projective plane) in terms of 35 forbidden minors~\cite{archdeacon1981kuratowski};
	\item characterisation of outer projective planar graphs (graphs that are embeddable on the projective plane with a disk removed such that all the vertices are located on the boundary of the surface) in terms of 32 forbidden minors~\cite{archdeacon1998obstruction};
	\item characterisation of outercylindrical graphs (graphs that are embeddable on a plane with two disks removed from it such that all the vertices are located on the boundary of the removed disks) in terms of 38 forbidden minors~\cite{archdeacon2001obstruction}.
\end{itemize}

%Graph $G$ is a non-separating planar graph if there is a planar drawing $D$ of $G$ on the plane such that for any cycle $C$ and any two vertices that are not a vertex of $C$ in $G$, the two vertices are on the same side of $C$ in $D$.

%%%%%%%%%%%%%%%%%%%%%%%%%%%%%%%%%%%%%%%%
%
%  Preliminary Lemmas
%
%%%%%%%%%%%%%%%%%%%%%%%%%%%%%%%%%%%%%%%%

\section{Preliminary Lemmas}
\label{se:lemmas}

A path $P$ in a graph $G$ is said to be \emph{chordless} if there is no edge between any two non-consecutive vertices of $P$ in $G$. A \emph{$uv$-path} is a path from a vertex $u$ to a vertex $v$.

Vertices $u$ and $v$, in a subdivision $S$ of $K_{2,3}$, are called the \emph{terminal vertices} of $S$ if both $u$ and $v$ have degree 3 in $S$. Define the \emph{terminal paths} in $S$ as the three $uv$-paths in $S$.

\begin{figure}[tbh]
  \centering
    \includegraphics[width=0.2\textwidth]{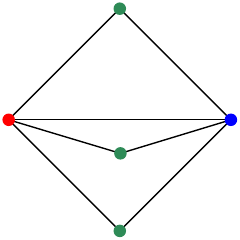}
  \caption{$K_{1,1,3}$}
  \label{fig:K_113}
\end{figure}

Next we will prove a couple of lemmas about the graphs that do not contain $K_{1,1,3}$ as a minor (see Figure~\ref{fig:K_113}). 

\begin{lemma}
\label{le:helper_lemma_1_for_characterisation_of_non-separating_planar_graphs}
	Every terminal path in a spanning $K_{2,3}$-subdivision of a $K_{1,1,3}$-minor-free graph is chordless.
\end{lemma}

\begin{proof}
	Suppose that such a terminal path $P$ has a chord $e$. Then it is easy to find a $K_{1,1,3}$ minor in the graph.
\end{proof}

A vertex $w$ of a $uv$-path $P$ is an \emph{inner vertex} of $P$ if $w \neq u$ and $w \neq v$. An edge $e$ of a path $P$ is an \emph{inner edge} of $P$ if $e$ is incident with two inner vertices of $P$. 
	
\begin{figure}[tbh]
  \centering
    \includegraphics[width=0.28\textwidth]{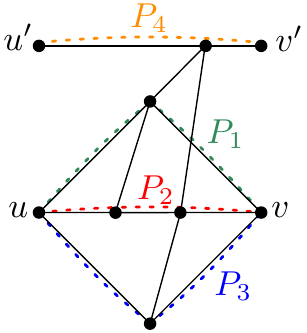}
  \caption{$P_2$ is the only middle path among the four paths $P_1, P_2, P_3, P_4$, where $P_1, P_2, P_3$ are $uv$-paths and $P_4$ is a $u'v'$-path.}
  \label{fig:middle_path}
\end{figure}

Given a set $\mathcal{P}$ of paths in a graph $G$, define a \emph{middle path} $P \in \mathcal{P}$ to be a path such that for any other path $P' \in \mathcal{P}$ there is an edge in $G$ that is incident with an inner vertex of $P$ and an inner vertex of $P'$. In other words, for each path $P' \in \mathcal{P}$ other than $P$ there is an inner vertex of $P$ that is adjacent to an inner vertex of $P'$ (see, e.g., Figure~\ref{fig:middle_path}). Two vertices $u$ and $v$ are \emph{co-path} with respect to $\mathcal{P}$ if $u$ and $v$ are on the same path in $\mathcal{P}$.

Any graph $G$ that contains a $K_{2,3}$-subdivision is \emph{middle-less} if there is no middle path among the terminal paths of any of the spanning $K_{2,3}$-subdivisions in $G$. Any graph $G$ with a spanning $K_{2,3}$-subdivision is \emph{middle-ful} if it is not middle-less.

We divide the rest of lemmas in this section into two subsections. The first section is about the middle-less graphs and the second section is about the middle-ful ones.

%%%%%%%%%%%%%%%%%%%%%%%%%%%%%%%%%%%%%%%%
%
%  Middle-less graphs
%
%%%%%%%%%%%%%%%%%%%%%%%%%%%%%%%%%%%%%%%%

\subsection{Middle-less Graphs}
\label{sse:middle-less_graphs}

We start by proving that middle-less graphs do not contain $W_4$ as a minor.

\begin{figure}[tbh]
  \centering
  \begin{subfigure}[t]{0.40\textwidth}
  \centering
    \includegraphics[width=0.5\textwidth]{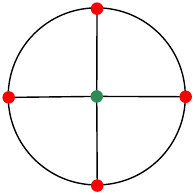}
    \caption{$W_4$}
    \label{sfig:W_4}
  \end{subfigure}
  \hspace{6 mm}
  \begin{subfigure}[t]{0.40\textwidth}
  \centering
    \includegraphics[width=0.5\textwidth]{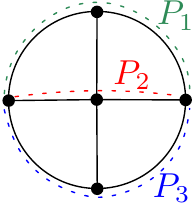}
    \caption{$W_4$ with $P_2$ as its middle path.}
    \label{sfig:W_4_and_middle_path}
  \end{subfigure}
  \caption{Any graph with a $W_4$ minor is middle-ful.}
  \label{fig:helper_lemma_6_for_characterisation_of_non-separating_planar_graphs}
\end{figure}

\begin{lemma}
	\label{le:helper_lemma_6_for_characterisation_of_non-separating_planar_graphs}
	If $G$ is a middle-less graph then $G$ does not contain $W_4$ as a minor (see Figure~\ref{sfig:W_4}).
\end{lemma}

\begin{proof}
	Suppose that there is a middle-less graph $G$ that contains $W_4$ as a minor. Then it is straightforward to find a $K_{2,3}$-subdivision with a middle path in $G$. But this is a contradiction since $G$ is middle-less (see, e.g., Figure~\ref{sfig:W_4_and_middle_path}).
\end{proof}

Let $U$ be a subset of the vertices of a graph $G$, then $G[U]$ denotes the subgraph of $G$ induced by $U$. Similarly, for any subgraph $H$ of the graph $G$, $G[H]$ denotes the subgraph of $G$ that is induced by the vertices of $H$.

\begin{lemma}
	\label{le:helper_lemma_5_for_characterisation_of_non-separating_planar_graphs}
	Let $P_1, P_2, P_3$ be the terminal paths in a spanning $K_{2,3}$-subdivision $S$ of a middle-less graph $G$ with no $K_{1,1,3}$-minor or $(K_1 \cup K_{2,3})$-minor where $G[P_1 \cup P_2]$ has an edge $e$ that is not in $P_1$ or $P_2$.
	Then:
	\begin{itemize}
%		\item One of $G[P_2 \cup P_3]$ and $G[P_3 \cup P_1]$ does not have any edge that is not in $P_1,P_2$ or $P_3$ and
		\item either every edge of $G[P_2 \cup P_3]$ is an edge of $P_2 \cup P_3$ or every edge of $G[P_2 \cup P_1]$ is an edge of $P_1 \cup P_3$ and
		\item $e$ is the only edge in $G[P_1 \cup P_2]$ that is not in $P_1, P_2$ and $P_3$. 
	\end{itemize}
\end{lemma}

\begin{proof}
	Let $G_1 = G[P_1 \cup P_2]$, $G_2 = G[P_2 \cup P_3]$ and $G_3 = G[P_3 \cup P_1]$ and let $u$ and $v$ be the two vertices of $e$. First we show that $G_2$ does not have any edge that is not an edge of $P_2$ or $P_3$. To reach a contradiction suppose that $G_2$ has an edge $e_1=(u_1,v_1)$ that is not in $P_2 \cup P_3$. Moreover, by the assumptions of the lemma, there is an edge $e$ in $G_1$ that is not in $P_1 \cup P_2$. 
	
	By Lemma~\ref{le:helper_lemma_1_for_characterisation_of_non-separating_planar_graphs}, $e$ and $e_1$ are not chords of $P_1,P_2$ or $P_3$ and therefore, without loss of generality, $u$ is an inner vertex of $P_1$ and $v$ is an inner vertex of $P_2$ and $u_1$ is an inner vertex of $P_2$ and $v_1$ is an inner vertex of $P_3$ (see, e.g., Figure~\ref{fig:helper_lemma_5_for_characterisation_of_non-separating_planar_graphs_1}). But this is a contradiction since then $P_2$ is a middle path and therefore $G$ is not middle-less. Similarly we can show that $G_3$ does not have any edge that is not an edge of $P_1$ or $P_3$. 
	
\begin{figure}[tbh]
  \centering
    \includegraphics[width=0.26\textwidth]{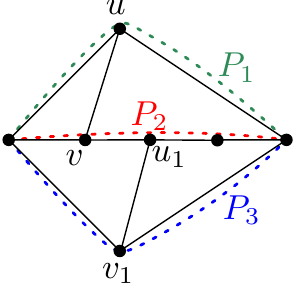}
  \caption{$u,v,u_1,v_1$ in $G$}
  \label{fig:helper_lemma_5_for_characterisation_of_non-separating_planar_graphs_1}
\end{figure}
	
	Now we show that there is at most one edge in $G_1$ that is not an edge of $P_1$ or $P_2$. To reach a contradiction suppose that $G_1$ has two edges $e_1=(u_1,v_1)$ and $e_2=(u_2,v_2)$ that are not among the edges of $P_1$ or $P_2$ (note that it is possible that either $u_1= u_2$ or $v_1= v_2$ or $u_1= v_2$ or $v_1= u_2$).
	
	By Lemma~\ref{le:helper_lemma_1_for_characterisation_of_non-separating_planar_graphs}, $e_1$ and $e_2$ are not chords of $P_1$ or $P_2$ and therefore, without loss of generality, let $u_1$ and $u_2$ be among the inner vertices of $P_1$ and $v_1$ and $v_2$ be among the inner vertices of $P_2$ (see, e.g., Figure~\ref{fig:helper_lemma_5_for_characterisation_of_non-separating_planar_graphs_2}). 
	
\begin{figure}[tbh]
  \centering
    \includegraphics[width=0.26\textwidth]{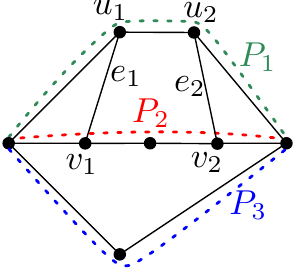}
  \caption{$e_1,e_2,P_1,P_2$ and $P_3$ in $G$}
  \label{fig:helper_lemma_5_for_characterisation_of_non-separating_planar_graphs_2}
\end{figure}
	
	Choose $P$ to be either $P_1$ or $P_2$ so that the endpoints of $e_1$ and $e_2$ on the other path are distinct. Let $G^-$ be the graph that is obtained by contracting all the edges of $P$ except the ones that are incident to the terminal vertices of $S$ into a single vertex $w$. It is easy to see that there is a $W_4$-minor in $G^-$ (see, e.g., Figure~\ref{fig:helper_lemma_5_for_characterisation_of_non-separating_planar_graphs_3}). Then by Lemma~\ref{le:helper_lemma_6_for_characterisation_of_non-separating_planar_graphs}, $G$ is not middle-less, which is a contradiction.
	
\begin{figure}[tbh]
  \centering
  \begin{subfigure}[t]{0.36\textwidth}
  \centering
    \includegraphics[width=0.75\textwidth]{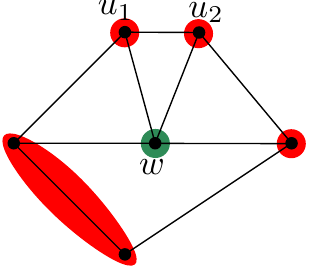}
    \caption{$G^-$}
    \label{sfig:helper_lemma_5_for_characterisation_of_non-separating_planar_graphs_4}
  \end{subfigure}
  \hspace{10 mm}
  \begin{subfigure}[t]{0.36\textwidth}
  \centering
    \includegraphics[width=0.75\textwidth]{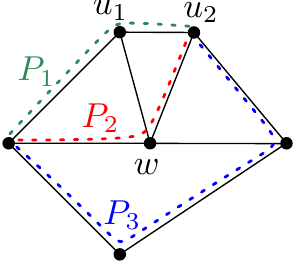}
    \caption{$G^-$ with a middle path $P_2$}
    \label{sfig:helper_lemma_5_for_characterisation_of_non-separating_planar_graphs_5}
  \end{subfigure}
  \caption{Finding a middle path in $G^-$.}
  \label{fig:helper_lemma_5_for_characterisation_of_non-separating_planar_graphs_3}
\end{figure}

\end{proof}

\begin{lemma}
	\label{le:helper_lemma_7_for_characterisation_of_non-sparse_non-separating_planar_graphs}
	Let $\{u, v\}$ and $\{P_1, P_2, P_3\}$ be the sets of terminal vertices and terminal paths respectively in a spanning $K_{2,3}$-subdivision $S$ of a middle-less graph $G$ with no $K_{1,1,3}$-minor and no $(K_1 \cup K_{2,3})$-minor where the lengths of $P_1$ and $P_2$ are greater than 2 and $G[P_1 \cup P_2]$ has an edge $e'=(u',v')$ that is not in $P_1 \cup P_2$. Then either:
	\begin{itemize}
		\item $u'$ and $v'$ are adjacent to $u$, or
		\item $u'$ and $v'$ are adjacent to $v$.
	\end{itemize}
\end{lemma}

\begin{proof}
	By Lemma~\ref{le:helper_lemma_1_for_characterisation_of_non-separating_planar_graphs}, $e'$ is not a chord of $P_1$ or $P_2$ and therefore, without loss of generality, let $u'$ be an inner vertex of $P_1$ and $v'$ be an inner vertex of $P_2$. To reach a contradiction, suppose that $u'$ and $v'$ are not both adjacent to the same vertex $u$ or $v$. We have two cases:
	
	\textbf{Case 1. Neither $u'$ nor $v'$ is adjacent to the terminal vertices.} In this case it is easy to find a $K_1 \cup K_{2,3}$ minor in $G$ (see, e.g., Figure~\ref{fig:helper_lemma_5_for_characterisation_of_non-separating_planar_graphs_6} and Figure~\ref{sfig:lemma_7_K_1_and_K_23_colour}).
	
\begin{figure}[tbh]
  \centering
  \begin{subfigure}[t]{0.36\textwidth}
  \centering
    \includegraphics[width=0.75\textwidth]{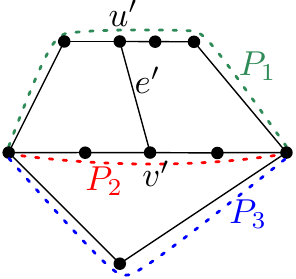}
    \caption{$e', P_1, P_2$ and $P_3$ in $G$}
    \label{sfig:helper_lemma_5_for_characterisation_of_non-separating_planar_graphs_9}
  \end{subfigure}
  \hspace{1 mm}
  \begin{subfigure}[t]{0.28\textwidth}
  \centering
    \includegraphics[width=0.99\textwidth]{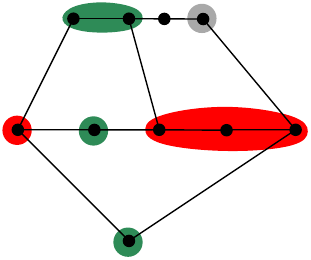}
    \caption{$P'_1$ and $P''_1$ in $G$}
    \label{sfig:helper_lemma_5_for_characterisation_of_non-separating_planar_graphs_10}
  \end{subfigure}
  \hspace{1 mm}
  \begin{subfigure}[t]{0.33\textwidth}
  \centering
    \includegraphics[width=0.90\textwidth]{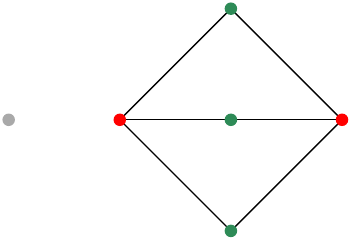}
    \caption{$K_1 \cup K_{2,3}$}
    \label{sfig:lemma_7_K_1_and_K_23_colour_1}
  \end{subfigure}
  \caption{$e', P_1, P_2, P_3, P'_1$ and $P''_1$ in $G$. Compare the colouring scheme of Figure~\ref{sfig:helper_lemma_5_for_characterisation_of_non-separating_planar_graphs_10} with Figure~\ref{sfig:lemma_7_K_1_and_K_23_colour} to see how $K_1 \cup K_{2,3}$ is a minor of $G$.}
  \label{fig:helper_lemma_5_for_characterisation_of_non-separating_planar_graphs_6}
\end{figure}

	\textbf{Case 2. One of the two vertices $u'$ or $v'$ is adjacent to $u$ or $v$.} Without loss of generality let $u'$ be adjacent to $u$ (see, e.g., Figure~\ref{sfig:helper_lemma_5_for_characterisation_of_non-separating_planar_graphs_6}).

\begin{figure}[tbh]
  \centering
  \begin{subfigure}[t]{0.36\textwidth}
  \centering
    \includegraphics[width=0.75\textwidth]{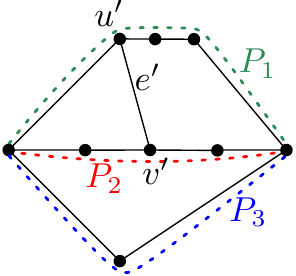}
    \caption{$e', P_1, P_2$ and $P_3$ in $G$}
    \label{sfig:helper_lemma_5_for_characterisation_of_non-separating_planar_graphs_6}
  \end{subfigure}
  \hspace{10 mm}
  \begin{subfigure}[t]{0.36\textwidth}
  \centering
    \includegraphics[width=0.82\textwidth]{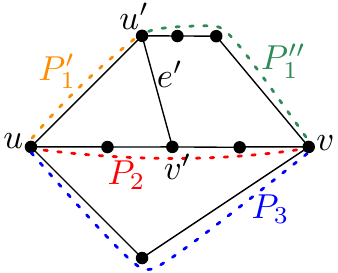}
    \caption{$P'_1$ and $P''_1$ in $G$}
    \label{sfig:helper_lemma_5_for_characterisation_of_non-separating_planar_graphs_7}
  \end{subfigure}
  \caption{$e', P_1, P_2, P_3, P'_1$ and $P''_1$ in $G$.}
  \label{fig:helper_lemma_5_for_characterisation_of_non-separating_planar_graphs_4}
\end{figure}
	The vertex $u'$ splits $P_1$ into two shorter paths $P'_1$ and $P''_1$, where $P_1$ consists of the edge $(u,u')$. Without loss of generality, let $P'_1$ be a shortest path among $P'_1$ and $P''_1$ (see, e.g., Figure~\ref{sfig:helper_lemma_5_for_characterisation_of_non-separating_planar_graphs_7}). Then, since  the lengths of $P_1$ and $P_2$ are greater than 2, it is easy to see that there is a $K_{2,3}$ minor in $P'_1 \cup e' \cup P_2 \cup P_3$ and an inner vertex $v''$ on $P''_1$ such that $P'_1 \cup e' \cup P_2 \cup P_3$ and $v''$ form a $K_1 \cup K_{2,3}$ minor in $G$ (see, e.g., Figure~\ref{fig:helper_lemma_5_for_characterisation_of_non-separating_planar_graphs_5}). However, this is a contradiction since $G$ is a $K_1 \cup K_{2,3}$-minor free graph.
	
\begin{figure}[tbh]
  \centering
    \includegraphics[width=0.30\textwidth]{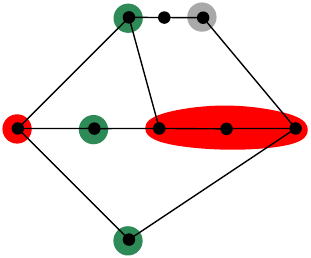}
  \caption{Finding a $K_1 \cup K_{2,3}$-minor in $G$ (compare with Figure~\ref{sfig:lemma_7_K_1_and_K_23_colour_1}).}
  \label{fig:helper_lemma_5_for_characterisation_of_non-separating_planar_graphs_5}
\end{figure}

\end{proof}

\begin{figure}[tbh]
  \centering
  \begin{subfigure}[t]{0.32\textwidth}
  \centering
    \includegraphics[width=0.80\textwidth]{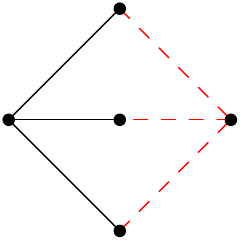}
    \caption{Type I}
    \label{sfig:characterisation_of_sparse_non-separating_planar_graphs_1}
  \end{subfigure}
  \hspace{0 mm}
  \begin{subfigure}[t]{0.32\textwidth}
  \centering
    \includegraphics[width=0.80\textwidth]{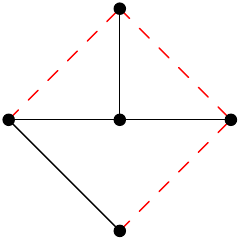}
    \caption{Type II}
    \label{sfig:characterisation_of_sparse_non-separating_planar_graphs_2}
  \end{subfigure}
  \hspace{0 mm}
  \begin{subfigure}[t]{0.32\textwidth}
  \centering
    \includegraphics[width=0.80\textwidth]{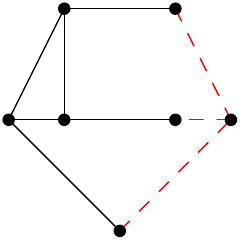}
    \caption{Type III}
    \label{sfig:characterisation_of_sparse_non-separating_planar_graphs_3}
  \end{subfigure}
  \caption{Three types of middle-less non-separating planar graphs}
  \label{fig:characterisation_of_sparse_non-separating_planar_graphs}
\end{figure}

\begin{lemma}
	\label{le:characterisation_of_sparse_non-separating_planar_graphs}
	Let $\mathcal{G}$ be the family of middle-less graphs with no $K_{1,1,3}$-minor, no $(K_1 \cup K_4)$-minor, no $(K_1 \cup K_{2,3})$-minor, and that contain a $K_{2,3}$-subdivision. Then any $G \in \mathcal{G}$ can be obtained by subdividing the red dashed edges of the graphs that are shown in Figure~\ref{fig:characterisation_of_sparse_non-separating_planar_graphs}.
\end{lemma}

\begin{proof}	
	Let $P_1,P_2,P_3$ be the terminal paths and $u,v$ be the terminal vertices in a $K_{2,3}$-subdivision $S$ of a graph $G \in \mathcal{G}$. Since $G$ does not contain $K_1 \cup K_{2,3}$ as a minor, $S$ is a spanning $K_{2,3}$-subdivision of $G$. If $G$ does not have any edges other than the edges of $P_1,P_2,P_3$ then, clearly, $G$ can be obtained by subdividing the red dashed edges of the graph depicted in Figure~\ref{sfig:characterisation_of_sparse_non-separating_planar_graphs_1}.
	
	Now let us consider the case where $G$ has an edge $e'=(u',v')$ that is not an edge of any of $P_1,P_2,P_3$. By Lemma~\ref{le:helper_lemma_5_for_characterisation_of_non-separating_planar_graphs}, $e'$ is the only edge in $G$ that is not an edge of $P_1,P_2$ or $P_3$. By Lemma~\ref{le:helper_lemma_1_for_characterisation_of_non-separating_planar_graphs}, $e'$ is not a chord of $P_1$ or $P_2$ and therefore, without loss of generality, let $u'$ be an inner vertex of $P_1$ and $v'$ be an inner vertex of $P_2$. We have two cases:
	
	\textbf{Case 1. Either $P_1$ or $P_2$ has length 2.} It is easy to verify that in this case $G$ is a graph that can be obtained by subdividing the red dashed edges in Figure~\ref{sfig:characterisation_of_sparse_non-separating_planar_graphs_2}.
	
	\textbf{Case 2. The lengths of both $P_1$ and $P_2$ are more than 2.} By Lemma~\ref{le:helper_lemma_7_for_characterisation_of_non-sparse_non-separating_planar_graphs}, both $u'$ and $v'$ are adjacent to the same vertex $u$ or $v$. Now it is easy to verify that in this case $G$ is a graph that can be obtained by subdividing the red dashed edges in Figure~\ref{sfig:characterisation_of_sparse_non-separating_planar_graphs_3}.
\end{proof}

%%%%%%%%%%%%%%%%%%%%%%%%%%%%%%%%%%%%%%%%
%
%  Middle-ful graphs
%
%%%%%%%%%%%%%%%%%%%%%%%%%%%%%%%%%%%%%%%%

\subsection{Middle-ful Graphs}
\label{sse:middle_full_graphs}

\begin{lemma}
\label{le:helper_lemma_2_for_characterisation_of_non-separating_planar_graphs}
	There is at most one middle path in the set of terminal paths of a spanning $K_{2,3}$-subdivision of a $K_{1,1,3}$-minor-free graph.
\end{lemma}

\begin{figure}[tbh]
  \centering
  \begin{subfigure}[t]{0.38\textwidth}
  \centering
    \includegraphics[width=0.65\textwidth]{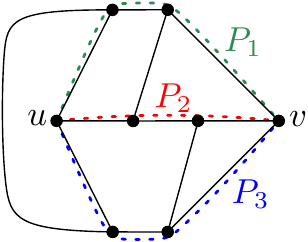}
    \caption{$G$ with $P_1$ and $P_2$ as middle paths.}
    \label{sfig:helper_lemma_2_for_characterisation_of_non-separating_planar_graphs_1}
  \end{subfigure}
  \hspace{10 mm}
  \begin{subfigure}[t]{0.38\textwidth}
  \centering
    \includegraphics[width=0.65\textwidth]{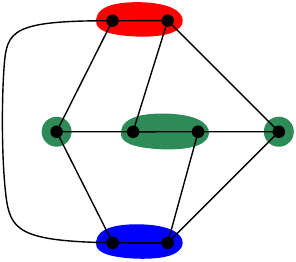}
    \caption{$G$ contains $K_{1,1,3}$ as a minor.}
    \label{sfig:helper_lemma_2_for_characterisation_of_non-separating_planar_graphs_2}
  \end{subfigure}
  \caption{If $P_1$ and $P_2$ are middle paths then $G$ contains $K_{1,1,3}$ as a minor. The colour scheme used here to colour the vertices of a $K_{1,1,3}$ minor is the same as the one used in Figure~\ref{fig:K_113}}
  \label{fig:helper_lemma_2_for_characterisation_of_non-separating_planar_graphs}
\end{figure}

\begin{proof}
	Let $\mathcal{P} = \{ P_1, P_2, P_3 \}$ be the set of terminal paths in a spanning $K_{2,3}$-subdivision $S$ in a $K_{1,1,3}$-minor-free graph $G$. To reach a contradiction, suppose that there is more than one middle path in $\mathcal{P}$. Without loss of generality, let $P_1$ and $P_2$ both be middle paths. Since $P_1$ and $P_2$ are middle paths:
	\begin{enumerate}
		\item there is an edge incident with an inner vertex of $P_1$ and an inner vertex of $P_2$, and
		\item there is an edge incident with an inner vertex of $P_1$ and an inner vertex of $P_3$, and 
		\item there is an edge incident with an inner vertex of $P_2$ and an inner vertex of $P_3$.
	\end{enumerate}
	Now, it is easy to find a $K_{1,1,3}$ as a minor in $G$. see, e.g., Figure~\ref{fig:helper_lemma_2_for_characterisation_of_non-separating_planar_graphs}.
\end{proof}

\begin{figure}[tbh]
  \centering
    \includegraphics[width=0.26\textwidth]{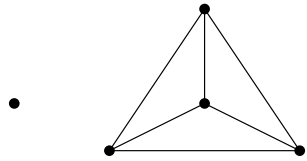}
  \caption{$K_1 \cup K_4$}
  \label{fig:K_1_and_K_4_colour}
\end{figure}

Next we will prove a lemma about a class of graphs that does not contain $K_1 \cup K_4$ as a minor (see Figure~\ref{fig:K_1_and_K_4_colour}).

\begin{lemma}
	\label{le:helper_lemma_3_for_characterisation_of_non-separating_planar_graphs}
	Let $P_1, P_2, P_3$ be the terminal paths in a spanning $K_{2,3}$-subdivision $S$ of a graph $G$ with no $(K_1 \cup K_4)$-minor, where $P_2$ is a middle path. Then there is no pair of edges $e_1=(u_1,v_1)$ and $e_2=(u_1,v_2)$ in $G$ such that $u_1$ is an inner vertex of $P_1$ or $P_3$ and $v_1$ and $v_2$ are two distinct inner vertices of $P_2$.
\end{lemma}

\begin{proof}
	To reach a contradiction suppose that there is an edge $e_1=(u_1,v_1)$ and an edge $e_2=(u_1,v_2)$ such that $u_1$ is an inner vertex of $P_1$ or $P_3$ and $v_1$ and $v_2$ are two inner vertices of $P_2$ (see, e.g., Figure~\ref{sfig:helper_lemma_3_for_characterisation_of_non-separating_planar_graphs_1}). Without loss of generality let $u_1$ be an inner vertex of $P_1$. Since $P_2$ is a middle path there is also an edge $e_3=(u_3,v_3)$ in $G$ such that $u_3$ is an inner vertex of $P_3$ and $v_3$ is an inner vertex of $P_2$ (see, e.g., Figure~\ref{sfig:helper_lemma_3_for_characterisation_of_non-separating_planar_graphs_2}). Now it is easy to find a $(K_1 \cup K_4)$-minor in $G$ (see, e.g., Figure~\ref{sfig:helper_lemma_3_for_characterisation_of_non-separating_planar_graphs_3}).
	
\begin{figure}[tbh]
  \centering
  \begin{subfigure}[t]{0.30\textwidth}
  \centering
    \includegraphics[width=0.97\textwidth]{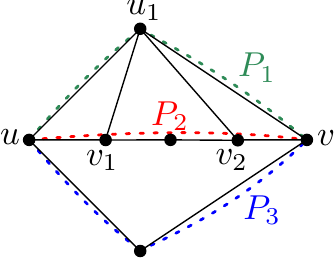}
    \caption{$G$ with $P_2$ as a middle path.}
    \label{sfig:helper_lemma_3_for_characterisation_of_non-separating_planar_graphs_1}
  \end{subfigure}
  \hspace{4 mm}
  \begin{subfigure}[t]{0.30\textwidth}
  \centering
    \includegraphics[width=0.97\textwidth]{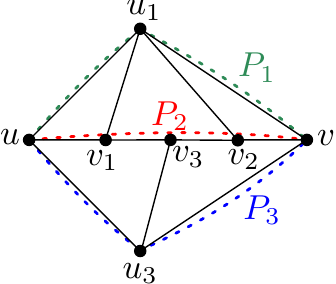}
    \caption{$G$ contains $K_{1,1,3}$ as a minor.}
    \label{sfig:helper_lemma_3_for_characterisation_of_non-separating_planar_graphs_2}
  \end{subfigure}
  \hspace{4 mm}
  \begin{subfigure}[t]{0.30\textwidth}
  \centering
    \includegraphics[width=0.97\textwidth]{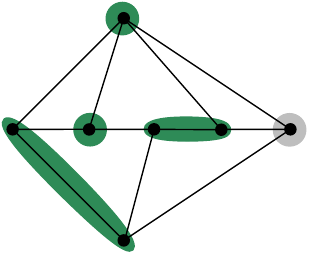}
    \caption{Finding a $K_1 \cup K_4$ minor in $G$.}
    \label{sfig:helper_lemma_3_for_characterisation_of_non-separating_planar_graphs_3}
  \end{subfigure}
  \caption{$G$ contains $K_1 \cup K_4$ as a minor.}
  \label{fig:helper_lemma_3_for_characterisation_of_non-separating_planar_graphs}
\end{figure}

\end{proof}

Let $P$ be a path and $h$ be a vertex that is not in $P$. Let $G$ be the graph that is obtained from $P$ and $h$ by adding an edge $(h,v)$ for every vertex $v$ in $P$. Then $G$ is a \emph{fan} graph and $h$ is the \emph{handle} of $G$.~\footnote{$K_3$ and $K_4$ minus an edge are the only fan graphs that do not have a unique handle.}

Let $P$ be a $uv$-path. We define the \emph{outer} inner vertices of $P$ as those inner vertices of $P$ that are adjacent to $u$ and $v$ on $P$.

\begin{lemma}
	\label{le:helper_lemma_7_for_characterisation_of_non-separating_planar_graphs}
	Let $P_1, P_2, P_3$ be the terminal paths in a spanning $K_{2,3}$-subdivision $S$ of a graph $G$ with no $K_{1,1,3}$-minor, no $(K_1 \cup K_4)$-minor and no $(K_1 \cup K_{2,3})$-minor, where $P_2$ is a middle path. Then, $G[P_1 \cup P_2]$ and $G[P_2 \cup P_3]$ are subgraphs of fan graphs whose handles are among the outer inner vertices of $P_2$.
%	\begin{enumerate}
%		\item $G[P_1 \cup P_2]$ and $G[P_2 \cup P_3]$ are subgraphs of fan graphs $G^+_1$ with handle $h_1$ and $G^+_2$ with handle $h_2$ and
%		\item $h_1$ and $h_2$ are the outer inner vertices of $P_2$.
%	\end{enumerate}
\end{lemma}

\begin{proof}
	Let $G_1 = G[P_1 \cup P_2]$ and $G_2 = G[P_2 \cup P_3]$. First we show that $G_1$ and $G_2$ are subgraphs of fan graphs. To reach a contradiction suppose that either $G_1$ or $G_2$ is not a subgraph of a fan graph. Without loss of generality, suppose that $G_1$ is not a subgraph of a fan graph. 
	
	Since $G_1$ is not a subgraph of a fan graph, there are two edges $e_1=(u_1,v_1)$ and $e_2=(u_2,v_2)$ in $G_1$ that are not an edge of $P_1$ or an edge of $P_2$ and are vertex-disjoint. By Lemma~\ref{le:helper_lemma_1_for_characterisation_of_non-separating_planar_graphs}, $e_1$ and $e_2$ are not chords of $P_1$ or $P_2$. In other words: 
	\begin{itemize}
		\item $u_1,v_1,u_2,v_2$ are all inner vertices of $P_1$ and $P_2$.
		\item $u_1$ and $v_1$ are not co-path with respect to $\{P_1, P_2, P_3\}$.
		\item $u_2$ and $v_2$ are not co-path with respect to $\{P_1, P_2, P_3\}$.
	\end{itemize}
	
	Without loss of generality let $u_1$ and $u_2$ be the two endpoints of $e_1$ and $e_2$ on $P_1$ and let $v_1$ and $v_2$ be the other two endpoints of $e_1$ and $e_2$ on $P_2$. Let $u$ and $v$ be the terminal vertices of $S$. Contract all the edges of $P_1$ that are not incident to $u$ and $v$ into a single vertex $w$ and let us denote the resulting minor of $G$ by $H$. 
	
	Since $H$ is a minor of $G$, it does not contain a $K_1 \cup K_4$ minor. Moreover, $P_2$ is a middle path in $H$. Also, $w$ is adjacent to $v_1$ and $v_2$ in $H$. Therefore, $e_1=(w,v_1)$ and $e_2=(w,v_2)$ are two edges of $H$ that contradict Lemma~\ref{le:helper_lemma_3_for_characterisation_of_non-separating_planar_graphs} and therefore $G_1$ is a subgraph of a fan graph. We denote the corresponding fan graph by $G^+_1$.
	
	Similarly, we conclude that $G_2$ is a subgraph of a fan graph and we denote the corresponding fan graph by $G^+_2$.
	
	Next we show that the handles of fan graphs $G^+_1$ and $G^+_2$, which we denote by $h_1$ and $h_2$ respectively, are outer inner vertices of $P_2$. As the first step, we show that $h_1$ and $h_2$ are inner vertices of $P_2$ and then as the second step we show that both $h_1$ and $h_2$ are adjacent to either $u$ or $v$ on $P_2$ (i.e., $h_1$ and $h_2$ are outer inner vertices of $P_2$).
	
	We use contradiction to prove the first step. To reach a contradiction suppose that either the handle of $G^+_1$ or the handle of $G^+_2$ is not an inner vertex of $P_2$. Without loss of generality, suppose that the handle of $G^+_1$ is not an inner vertex of $P_2$. Then it must be on $P_1$. So there are two edges $e_1=(u'_1,v')$ and $e_2=(u'_2,v')$ in $G_1$ that are not in $E(P_1) \cup E(P_2)$ and are incident with the same vertex $v'$ on $P_1$.
	
	By Lemma~\ref{le:helper_lemma_1_for_characterisation_of_non-separating_planar_graphs}, $e_1$ and $e_2$ are not chords of $P_1$ or $P_2$ and therefore $v'$ is an inner vertex of $P_1$ and $u'_1$ and $u'_2$ are inner vertices of $P_2$. However, this is also in contradiction with Lemma~\ref{le:helper_lemma_3_for_characterisation_of_non-separating_planar_graphs}.
	
	We use contradiction to prove the second step as well. To reach a contradiction, without loss of generality, suppose that $h_1$ is not adjacent to $u$ or $v$ on $P_2$ and let $h_2$ be any vertex on $P_2$. The handle $h_2$ splits $P_2$ into two subpaths: $P'_2$ from $u$ to $h_2$ and $P''_2$ from $h_2$ to $v$. Without loss of generality, let $h_1$ be on $P'_2$ or let $h_1 = h_2$ (see, e.g., Figure~\ref{sfig:helper_lemma_7_for_characterisation_of_non-separating_planar_graphs_1}).
	
	Since $P_2$ is a middle path, there are two edges $e_1= (u_1,x_1)$ and $e_2= (u_2,x_2)$ such that $u_1$ is an inner vertex of $P_1$ and $u_2$ is an inner vertex of $P_3$. Since $G_1$ is a subgraph of a fan graph $G^+_1$ with handle $h_1$ we have $x_1 = h_1$ and since $G_2$ is a subgraph of fan graph $G^+_2$ with handle $h_2$ we have $x_2 = h_2$ (see, e.g., Figure~\ref{sfig:helper_lemma_7_for_characterisation_of_non-separating_planar_graphs_2}). Let $P'_1$ be the part of $P_1$ from $u$ to $u_1$ and let $P'_3$ be the part of $P_3$ from $u$ to $u_2$. 
	
%	Since $P_2$ is a middle path and $G_1$ is a subgraph of a fan graph $G^+_1$ with handle $h_1$ and $G_2$ is a subgraph of fan graph $G^+_2$ with handle $h_2$, there are two edges $e_1= (u_1,h_1)$ and $e_2= (u_2,h_2)$ such that $u_1$ is an inner vertex of $P_1$ and $u_2$ is an inner vertex of $P_3$ (see, e.g., Figure~\ref{sfig:helper_lemma_7_for_characterisation_of_non-separating_planar_graphs_2}). Let $P'_1$ be the part of $P_1$ from $u$ to $u_1$ and let $P'_3$ be the part of $P_3$ from $u$ to $u_2$. 

\begin{figure}[tbh]
  \centering
  \begin{subfigure}[t]{0.35\textwidth}
  \centering
    \includegraphics[width=0.97\textwidth]{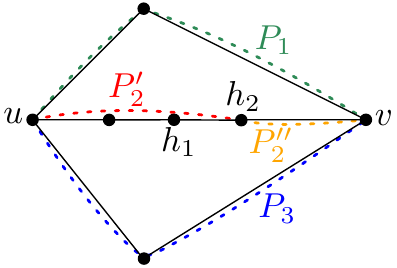}
    \caption{$h_1,h_2,P'_2$ and $P''_2$ in $G$.}
    \label{sfig:helper_lemma_7_for_characterisation_of_non-separating_planar_graphs_1}
  \end{subfigure}
  \hspace{8 mm}
  \begin{subfigure}[t]{0.35\textwidth}
  \centering
    \includegraphics[width=0.97\textwidth]{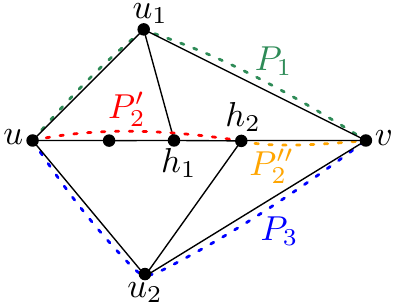}
    \caption{$u_1,u_2$ in $G$.}
    \label{sfig:helper_lemma_7_for_characterisation_of_non-separating_planar_graphs_2}
  \end{subfigure}
  \caption{Finding $K_1 \cup K_{2,3}$ minor in $G$.}
  \label{fig:helper_lemma_7_for_characterisation_of_non-separating_planar_graphs_1}
\end{figure}
	
	Now it is easy to see that $v$ together with $P'_1 \cup (u_1,h_1) \cup P'_2 \cup (u_2,h_2) \cup P'_3$ contains a $K_1 \cup K_{2,3}$ minor, which is a contradiction (see, e.g., Figure~\ref{fig:helper_lemma_7_for_characterisation_of_non-separating_planar_graphs_2}).  
	
\begin{figure}[tbh]
  \centering
  \begin{subfigure}[t]{0.36\textwidth}
  \centering
    \includegraphics[width=0.97\textwidth]{figures/K_1_and_K_23_colour.pdf}
    \caption{$K_1 \cup K_{2,3}$}
    \label{sfig:lemma_7_K_1_and_K_23_colour}
  \end{subfigure}
  \hspace{8 mm}
  \begin{subfigure}[t]{0.35\textwidth}
  \centering
    \includegraphics[width=0.97\textwidth]{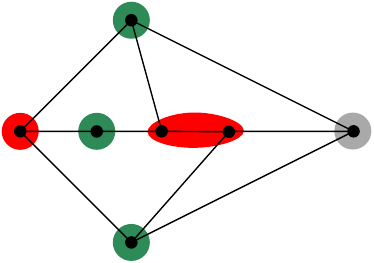}
    \caption{$G$ contains $K_1 \cup K_{2,3}$ as a minor.}
    \label{sfig:helper_lemma_7_for_characterisation_of_non-separating_planar_graphs_3}
  \end{subfigure}
  \caption{Finding $K_1 \cup K_{2,3}$ minor in $G$.}
  \label{fig:helper_lemma_7_for_characterisation_of_non-separating_planar_graphs_2}
\end{figure}

\end{proof}

\begin{lemma}
	\label{le:helper_lemma_8_for_characterisation_of_non-separating_planar_graphs}
	Let $G[P_1 \cup P_2]$ and $G[P_2 \cup P_3]$ be subgraphs of fan graphs $G^+_1$ and $G^+_2$ with the same handle $h$ where $P_1, P_2, P_3$ are the terminal paths in a spanning $K_{2,3}$-subdivision $S$ of a $K_{1,1,3}$-minor-free and $(K_1 \cup K_{2,3})$-minor-free graph $G$ in which $P_2$ is a middle path. Then length of $P_2$ is 2.
\end{lemma}

\begin{proof}
	To reach a contradiction suppose that length of $P_2$ is greater than 2. Since $P_2$ is the middle path, by Lemma~\ref{le:helper_lemma_7_for_characterisation_of_non-separating_planar_graphs}, $h$ is an outer inner vertex of $P_2$. Now, it is easy to find a $K_1 \cup K_{2,3}$ minor in $G$ which contradicts the assumptions of the lemma (see, e.g., Figure~\ref{fig:helper_lemma_8_for_characterisation_of_non-separating_planar_graphs_1}).
	
\begin{figure}[tbh]
  \centering
    \includegraphics[width=0.26\textwidth]{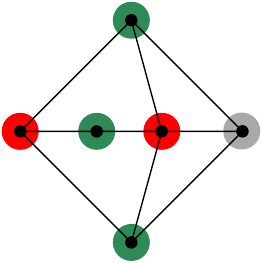}
  \caption{Finding $K_1 \cup K_{2,3}$ in $G$}
  \label{fig:helper_lemma_8_for_characterisation_of_non-separating_planar_graphs_1}
\end{figure}

\end{proof}

\begin{lemma}
	\label{le:helper_lemma_9_for_characterisation_of_non-separating_planar_graphs}
	Let $G_1 = G[P_1 \cup P_2]$ and $G_2 = G[P_2 \cup P_3]$ be subgraphs of fan graphs $G^+_1$ with handle $h_1$ and $G^+_2$ with handle $h_2$ respectively such that $h_1 \neq h_2$, where $u, v$ are the terminal vertices and $P_1, P_2, P_3$ are the terminal paths in a spanning $K_{2,3}$-subdivision $S$ of a graph $G$ with no $K_{1,1,3}$-minor, no $(K_1 \cup K_4)$-minor, no $(K_1 \cup K_{2,3})$-minor, and in which $P_2$ is a middle path. 
	
	Then there is exactly one edge $e'=(h_1,v')$ in $G_1$ that is not in $P_1 \cup P_2$ and there is exactly one edge $e''=(h_2,v'')$ in $G_2$ that is not in $P_2 \cup P_3$, where: 
	\begin{itemize}
		\item $h_1$ and $v'$ are outer inner vertices of $P_2$ and $P_1$ respectively that are both adjacent to $u$ or both adjacent to $v$ and
		\item $h_2$ and $v''$ are outer inner vertices of $P_2$ and $P_3$ respectively that are both adjacent to $u$ or both adjacent to $v$.
	\end{itemize}
\end{lemma}

\begin{proof}
	Since $P_2$ is a middle path, there is an edge $e'=(h_1,v')$ in $G_1$ that is not in $P_1 \cup P_2$ and there is an edge $e''=(h_2,v'')$ in $G_2$ that is not in $P_2 \cup P_3$. Moreover, by Lemma~\ref{le:helper_lemma_7_for_characterisation_of_non-separating_planar_graphs}, $h_1$ and $h_2$ are outer inner edges of $P_2$.

	Now to reach a contradiction, without loss of generality, let $h_1$ be adjacent to $u$ on $P_2$ but let $v'$ be a vertex that is not adjacent to $u$ on $P_1$. Let $v_1$ be the vertex that is adjacent to $u$ on $P_1$.
	
	Since $h_1$ and $h_2$ are inner vertices of the middle path $P_2$, by Lemma~\ref{le:helper_lemma_1_for_characterisation_of_non-separating_planar_graphs}, $v'$ is an inner vertex of $P_1$ and $v''$ is an inner vertex of $P_2$ (see, e.g., Figure~\ref{sfig:helper_lemma_9_for_characterisation_of_non-separating_planar_graphs_1}).
	
	We know that $v_1$ appears before $v'$ as we traverse $P_1$ from $u$ towards $v$ and $h_1$ appears before $h_2$ as we traverse $P_2$ from $u$ towards $v$. Let $P'$ be the part of $P_1$ that stretches from $v'$ to $v$. Now it is easy to see that $v_1$ together with $(h_1,v') \cup P' \cup P_2 \cup P_3$ contains a $K_1 \cup K_{2,3}$ minor, which is a contradiction (see, e.g., Figure~\ref{sfig:helper_lemma_9_for_characterisation_of_non-separating_planar_graphs_2}).  
	
\begin{figure}[tbh]
  \centering
  \begin{subfigure}[t]{0.30\textwidth}
  \centering
    \includegraphics[width=0.80\textwidth]{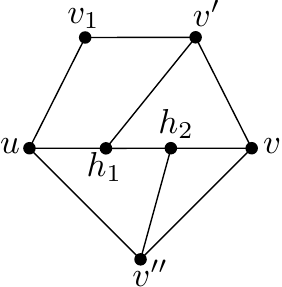}
    \caption{$K_1 \cup K_{2,3}$}
    \label{sfig:helper_lemma_9_for_characterisation_of_non-separating_planar_graphs_1}
  \end{subfigure}
  \hspace{8 mm}
  \begin{subfigure}[t]{0.30\textwidth}
  \centering
    \includegraphics[width=0.80\textwidth]{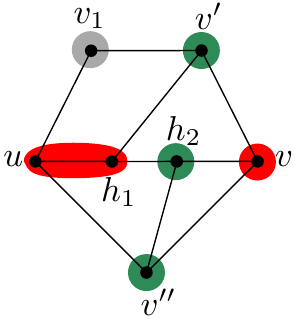}
    \caption{$G$ contains $K_1 \cup K_{2,3}$ as a minor.}
    \label{sfig:helper_lemma_9_for_characterisation_of_non-separating_planar_graphs_2}
  \end{subfigure}
  \caption{Finding a $K_1 \cup K_{2,3}$ minor in $G$.}
  \label{fig:helper_lemma_9_for_characterisation_of_non-separating_planar_graphs}
\end{figure}

\end{proof}

\begin{lemma}
	\label{le:characterisation_of_non-sparse_non-separating_planar_graphs}
	Let $\mathcal{G}$ be the family of middle-ful $K_{1,1,3}$-minor-free, $(K_1 \cup K_4)$-minor-free and $(K_1 \cup K_{2,3})$-minor-free graphs that contain a $K_{2,3}$-subdivision. Then any $G \in \mathcal{G}$ is either a subgraph of a wheel with at least 4 spokes or it is an elongated triangular prism.
\end{lemma}

\begin{proof}	
	Let $P_1,P_2,P_3$ be the terminal paths and $u,v$ be the terminal vertices in a $K_{2,3}$-subdivision $S$ of a graph $G \in \mathcal{G}$ where $P_2$ is a middle path. Since $G$ does not contain $K_1 \cup K_{2,3}$ as a minor, $S$ is a spanning $K_{2,3}$-subdivision of $G$. Let $G_1= G[P_1 \cup P_2]$ and $G_2= G[P_2 \cup P_3]$. Since $P_2$ is a middle path, by Lemma~\ref{le:helper_lemma_7_for_characterisation_of_non-separating_planar_graphs}, $G_1$ and $G_2$ are subgraphs of fan graph $G_1^+$ and $G_2^+$ with handles $h_1$ and $h_2$ where $h_1$ and $h_2$ are both among the outer inner vertices of $P_2$.
	
	We break the rest of the proof into two cases:
	
	\textbf{Case 1. $h_1 = h_2$.} By Lemma~\ref{le:helper_lemma_8_for_characterisation_of_non-separating_planar_graphs}, the length of $P_2$ is 2 and therefore $G$ is a subgraph of a wheel $W$. Moreover, since $P_2$ is a middle path, $W$ has at least 4 spokes. 
	
	\textbf{Case 2. $h_1 \neq h_2$.} By Lemma~\ref{le:helper_lemma_9_for_characterisation_of_non-separating_planar_graphs}, there is exactly one edge $e_1$ in $G_1$ that is not in $P_1 \cup P_2$ and exactly one edge $e_2$ in $G_2$ that is not in $P_2 \cup P_3$. Then $G$ is an elongated triangular prism. 
\end{proof}

%\begin{figure}[tbh]
%  \centering
%    \includegraphics[width=0.33\textwidth]{figures/prism.pdf}
%  \caption{A family of middle-ful non-separating planar graphs}
%  \label{fig:characterisation_of_non-sparse_non-separating_planar_graphs}
%\end{figure}

%%%%%%%%%%%%%%%%%%%%%%%%%%%%%%%%%%%%%%%%
%
%  Proof of the Main Theorems
%
%%%%%%%%%%%%%%%%%%%%%%%%%%%%%%%%%%%%%%%%

\section{Proof of the Main Theorems}
\label{se:proof}

\begin{lemma}
\label{le:construction_of_non-separating_planar_graphs}
	A graph $G$ does not contain any of $K_1 \cup K_4$ or $K_1 \cup K_{2,3}$ or $K_{1,1,3}$ as a minor if and only if $G$ is either an outerplanar graph or a subgraph of a wheel or an elongated triangular prism. 
\end{lemma}

\begin{proof}
	It is straightforward to see that any outerplanar graph or a subgraph of a wheel or an elongated triangular prism does not contain any of $K_1 \cup K_4$ or $K_1 \cup K_{2,3}$ or $K_{1,1,3}$ as a minor. Next we prove the lemma in the other direction.
	
	We break the proof into the following three cases:
	\begin{enumerate}
		\item $G$ does not contain any of $K_4$ or $K_{2,3}$ as a minor.
		\item $G$ contains $K_4$ but does not contain $K_{2,3}$ as a minor.
		\item $G$ contains $K_{2,3}$ as a minor.
	\end{enumerate}
	
	\textbf{Case 1. $G$ does not contain any of $K_4$ or $K_{2,3}$ as a minor.} In this case, $G$ is outerplanar.
	
	\textbf{Case 2. $G$ contains $K_4$ as a minor but it does not contain $K_{2,3}$ as a minor.} Since the degrees of the vertices in $K_4$ are less than 4, any subgraph contractible to $K_4$ is also a subdivision of $K_4$. Therefore, there is a subdivision $S$ of $K_4$ in $G$. 
	
	Since $G$ does not contain $K_1 \cup K_4$ as a minor, $S$ is a spanning subgraph of $G$ (any vertex of $G$ is also a vertex of $S$). Moreover, since any proper subdivision of $K_4$ contains $K_{2,3}$ as a minor, $K_4$ is the only graph that contains $K_4$ as a minor but does not contain $K_{2,3}$ as a minor. So $G$ is isomorphic to $K_4$ and is a subgraph of a wheel.
	
	\textbf{Case 3. $G$ contains $K_{2,3}$ as a minor.} Since the degrees of the vertices in $K_{2,3}$ are less than 4, any subgraph contractible to $K_{2,3}$ is also a subdivision of $K_{2,3}$. Therefore, there is a subdivision $S$ of $K_{2,3}$ in $G$. Since $G$ does not contain $K_1 \cup K_{2,3}$ as a minor, $S$ is a spanning subgraph of $G$. 	
	
	Here we have two cases:
	
	\textbf{Case 3a. Graph $G$ is middle-less.} By Lemma~\ref{le:characterisation_of_sparse_non-separating_planar_graphs}, $G$ can be obtained by subdividing the red dashed edges of one of the graphs shown in Figure~\ref{fig:characterisation_of_sparse_non-separating_planar_graphs}. Now, any of the graphs shown in Figure~\ref{fig:characterisation_of_sparse_non-separating_planar_graphs} is either a subgraph of a wheel or an elongated triangular prism. Therefore $G$ is either a subgraph of a wheel or a subgraph of an elongated triangular prism.
	
	\textbf{Case 3b. Graph $G$ is middle-ful.} By Lemma~\ref{le:characterisation_of_non-sparse_non-separating_planar_graphs}, $G$ is either a subgraph of a wheel or it is an elongated triangular prism.	
\end{proof}

%\begin{theorem}
%\label{th:characterisation_of_non-separating_planar_graphs}
%	$G$ is a non-separating planar graph if and only if $G$ does not contain any of $K_1 \cup K_4$ or $K_1 \cup K_{2,3}$ or $K_{1,1,3}$ as a minor.
%\end{theorem}

\subsection{Proof of Theorem~\ref{th:non-separating_planar_graphs_excluded_minors}}
Now we are ready to prove Theorem~\ref{th:non-separating_planar_graphs_excluded_minors}.

\begin{proof}
	It is straightforward to verify that in any planar drawing of a graph that contains $K_1 \cup K_4$ or $K_1 \cup K_{2,3}$ or $K_{1,1,3}$ as a minor, there are two vertices that are separated by a cycle. Therefore, to prove this theorem, it is sufficient to show that any graph that does not contain any of $K_1 \cup K_4$ or $K_1 \cup K_{2,3}$ or $K_{1,1,3}$ as a minor is a non-separating planar graph.
	
	By Lemma~\ref{le:construction_of_non-separating_planar_graphs}, any graph that does not contain any of $K_1 \cup K_4$ or $K_1 \cup K_{2,3}$ or $K_{1,1,3}$ as a minor is either an outerplanar graph or a subgraph of a wheel or an elongated triangular prism and it is easy to verify that any such graph is a non-separating planar graph.
\end{proof}

\subsection{Proof of Theorem~\ref{th:structural_characterisation_of_non-separating_planar_graphs}}

\begin{proof}
Theorem~\ref{th:structural_characterisation_of_non-separating_planar_graphs} is a direct consequence of Lemma~\ref{le:construction_of_non-separating_planar_graphs} and Theorem~\ref{th:non-separating_planar_graphs_excluded_minors}.
\end{proof}

Theorems~\ref{th:non-separating_planar_graphs_excluded_minors} and \ref{th:structural_characterisation_of_non-separating_planar_graphs} together provide us with Theorem~\ref{th:complete_characterisation_of_non-separating_planar_graphs}:

\begin{theorem}
\label{th:complete_characterisation_of_non-separating_planar_graphs}
	The following are equivalent, for any graph $G$:
	\begin{enumerate}
		\item $G$ does not contain any of $K_1 \cup K_4$ or $K_1 \cup K_{2,3}$ or $K_{1,1,3}$ as a minor.
		\item $G$ is outerplanar or a subgraph of a wheel or an elongated triangular prism.
		\item $G$ is a non-separating planar graph.
	\end{enumerate}
\end{theorem}

\section{Proof of Theorem~\ref{th:number_of_edges_in_flat_graphs}}
\label{se:linkless_edges_count}

\begin{figure}
  \centering
  \begin{subfigure}[t]{0.26\textwidth}
    \includegraphics[width=0.99\textwidth]{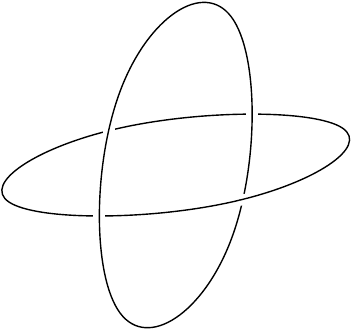}
    \caption{a link in 3D}
    \label{sfig:links_and_seprating_cycles_1}
  \end{subfigure}
  \begin{subfigure}[t]{0.47\textwidth}
    \includegraphics[width=0.99\textwidth]{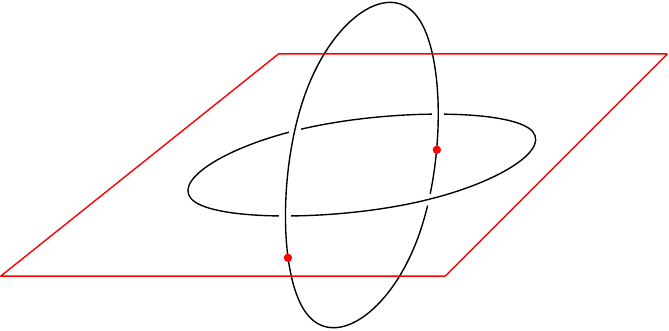}
    \caption{intersection of a plane with a link}
    \label{sfig:links_and_seprating_cycles_2}
  \end{subfigure}
  \hspace{2 mm}
  \begin{subfigure}[t]{0.22\textwidth}
    \includegraphics[width=0.9\textwidth]{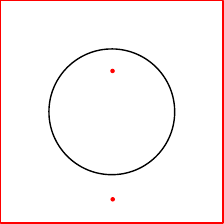}
    \caption{cross section of the link from the top}
    \label{sfig:links_and_seprating_cycles_3}
  \end{subfigure}
  \caption{A separating cycle in a cross section of a link with a plane}
  \label{fig:links_and_seprating_cycles}
\end{figure}

Consider two linked circles in three dimensions and a cross section of them that contains one of the two circles as depicted in Figure~\ref{sfig:links_and_seprating_cycles_1}. Such a cross-section has a structure that resembles the structure of a separating cycle with a vertex inside it and another outside it.

With this intuition in mind, we prove Theorem~\ref{th:number_of_edges_in_flat_graphs}.

\begin{proof}
	
Let $G$ be an elongated prism. $G$ is a maximal non-separating graph with $|V(G)| + 3$ edges. Moreover, $G$ contains both $K_4$ and $K_{2,3}$ as minors. 

Let $H$ be the graph that is obtained by adding two new vertices $u$ and $v$ to $G$ such that $u$ and $v$ are each adjacent to all the vertices of $G$. $H$ has at most $3|V(H)| - 3$ edges. We claim that $H$ is a maximal linkless graph. 

To prove that $H$ is a maximal linkless graph, we show that any graph $H^+$ that is obtained by adding an edge $e$ to $H$ is not a linkless graph. Since $u$ and $v$ are adjacent to all the vertices of $G$, edge $e$ (in $H^+$) is either $(u,v)$ or it is an edge between two vertices of $G$. 

Let $H^+$ be the graph obtained by adding $(u,v)$ to $H$. Since $G$ contains $K_4$ or $K_{2,3}$ as a minor, $H^+$ contains either $K_6$ or $K_{1,1,2,3}$ as a minor. The latter contains $K_{1,3,3}$. However, $K_6$ and $K_{1,3,3}$ are both forbidden minors for linkless graphs.

Now let $H^+$ be the graph that is obtained by adding an edge between two vertices of $G$ in $H$. By the characterisation of non-separating planar graphs, $G + e$ contains $K_4 \cup K_1$ or $K_{2,3} \cup K_1$ or $K_{1,1,3}$ as a minor. 

If $G + e$ contains $K_4 \cup K_1$ as a minor, then $H^+$ contains $K_6$ as a minor and hence $H^+$ is not a linkless graph.
If $G + e$ contains $K_{2,3} \cup K_1$ as a minor, then $H^+$ contains $K_{1,1,2,3}$ as a minor. But $K_{1,1,2,3}$ contains $K_{1,3,3}$ as a minor and therefore  $H^+$ is not a linkless graph.
If $G + e$ contains $K_{1,1,3}$ as a minor, then $H^+$ contains $K_{2,1,1,3}$ as a minor, which in turn contains $K_{1,3,3}$ as a minor. Therefore  $H^+$ is not a linkless graph.
\end{proof}

%%%%%%%%%%%%%%%%%%%%%%%%%%%%%%%%%%%%%%%%
%
%  Conclusion
%
%%%%%%%%%%%%%%%%%%%%%%%%%%%%%%%%%%%%%%%%

\section{Conclusion}
\label{se:conclusion}

This paper provides a forbidden minor characterisation for non-separating planar graphs and as expected, the forbidden minors for non-separating planar graphs are each a minor of one of the two forbidden minors for planar graphs.

Moreover it describes the structure of these graphs by proving that any maximal non-separating planar graph is either an outerplanar graph or a wheel or an elongated triangular prism.

One can define a similar class of graphs with respect to surfaces other than the plane. For example, a \emph{non-separating toroidal} graph is a graph that has a drawing $D$ on the torus such that: 
\begin{itemize}
	\item no two edges cross and 
	\item for any cycle $C$ in $D$ and any two vertices $u$ and $v$ in $D$ that are not a vertex of $C$, one can draw a curve from $u$ to $v$ without crossing any edge of $C$. %in the drawing that is obtained by restricting $D$ to the vertices and edges of $C \cup u \cup v$.
\end{itemize}
Any such class of graphs is also closed under minor operations and hence it can be characterised using a finite set of forbidden minors. It would be specially interesting to know the set of forbidden minors for non-separating toroidal graphs since they are all minors of the forbidden minors for toroidal graphs and we do not yet know the complete set of forbidden minors of toroidal graphs.

In Theorem~\ref{th:number_of_edges_in_flat_graphs}, we also showed that there are maximal linkless graphs with $3|V| -3$ edges. Now, a natural question that comes into mind is the following: Does every edge-maximal linkless graph have at least $3n-3$ edges?  Apart from that, Theorem~\ref{th:number_of_edges_in_flat_graphs} showed that there is a connection between non-separating planar graphs and linkless graphs. It would be interesting to explore this connection further. In fact it was this connection that served as our first motivation for exploring the structure of non-separating planar graphs.

Another application of non-separating planar graphs is in decomposing planar graphs. In another paper, we use such a decomposition to prove a stronger version of the Hanani-Tutte Theorem~\cite{deh19Str}. Finally, it would be interesting to see if there are other applications for non-separating planar graphs.

\section*{Acknowledgment}

The authors thank Peter Eades and David Wood for their helpful discussions and suggestions.

\bibliographystyle{abbrvurl}

%\bibliography{biblio}

\end{document}